# Convergence and Computation of a Class of Generalized Integrals

Haoding Meng

College of Life Science and Technology, Xi'an Jiaotong University, 710049, China

Email: menghd@stu.xjtu.edu.cn

April, 2019

**Abstract.** In this study, we discuss the convergence and divergence of generalized integrals $\int_{0}^{+\infty} \frac{\sin^{b} x}{x^{a}} dx \left( a \in R^{+}, b \in N^{+} \right)$, and use the transformation method, the partial integration method, the mathematical induction method and the complex variable residual method to derive the calculation formulas of the two types of integrals. The calculation formulas of the two types of integrals are unified. Finally, the convergence criterion and the unified calculation formula of $\int_{0}^{+\infty} \frac{\sin^{b} x}{x^{a}} dx$ is obtained.

## I  Introduction

Dirichlet integral, $\int_{0}^{+\infty} \frac{\sin x}{x} dx = \frac{\pi}{2}$, was given by Euler in 1781 through the residue method of complex variable function firstly. Robachevsky got a primary solution by using fractional decomposition. Chinese and foreign scholars have successively given a variety of ingenious solutions: double integral, abnormal integral with parameter variable etc[1].

There are endless researches on the generalization of this integral, among which Wu Zhiqin and other three got the calculation formula of $\int_{0}^{\infty} \left( \frac{\sin \alpha x}{x} \right)^{n} dx$ [2] for $n \geq 2, k \in N$, while Lin Feng and other three got the calculation formula of $\int_{-\infty}^{+\infty} \frac{\sin^{n} x}{x^{r}} dx$ [3] when $n$ and $r$ are both even or odd with $n \geq r > 0$.

Most of the previous studies have only discussed the corresponding calculation formula of $I = \int_0^{+\infty} \frac{\sin^b x}{x^a} dx$ when the power of the numerator denominator is positive integer and has the same parity, but not enough attention has been paid to its general power, and the convergence and divergence of this kind of integral are mostly ignored or some restrictions are given without being discussed.

In this paper, we firstly give the convergence results of the generalized integral $I = \int_0^{+\infty} \frac{\sin^b x}{x^a} dx$ for $a \in R^+, b \in N^+$, and give the calculation formula when $a$ is a positive integer. Then, we give the general calculation formula of $I$ by the partial integral method, trigonometric polynomial, Euler partial summation formula. It is a concise provision: the lower case letters in the following are all non negative integers, and $\Gamma(\cdot)$ means Gamma function.

## II  Convergence of generalized integral

**Lemma 1.** When $a \geq b+1 \geq 1$, $I = \int_0^{+\infty} \frac{\sin^b x}{x^a} dx$ diverges.

Proof.

Because

$$I = \int_0^{+\infty} \frac{\sin^b x}{x^a} dx = \int_0^1 \frac{\sin^b x}{x^a} dx + \int_1^{+\infty} \frac{\sin^b x}{x^a} dx,$$

and the first integral on the right side of the equation

$$\frac{\sin^b x}{x^a} \sim \frac{1}{x^{a-b}} \quad (x \to 0^+).$$

So for $a \geq b+1 \geq 1$, this part diverges and $I$ also diverges.

The following general always assumes $0 < a < b+1$. In order to obtain the convergence and divergence of $I$, it is necessary to use the following equations[4]

$$\sin^{2n-1} x = \frac{1}{2^{2n-2}} \sum_{k=0}^{n-1} (-1)^{n-k} C_{2n-1}^k \sin(2n-2k-1)x$$

$$\sin^{2n} x = \frac{1}{2^{2n}} \left\{ \sum_{k=0}^{n-1} (-1)^{n-k} 2C_{2n}^k \cos 2(n-k)x + C_{2n}^n \right\}.$$

And for any positive real number $p$, $\int_0^{+\infty} \frac{\sin px}{x} dx = \frac{\pi}{2}$.

If real number $a, b$ satisfies $0 < a < b+1$, the following conclusion is true.

(i) $I = \int_0^{+\infty} \frac{\sin^b x}{x^a} dx$ diverges when $a > 1$.

(ii) When $b$ is integer and $0 < a \le 1$, if $b$ is an odd number, $I$ converges, else if $b$ is an even number, $I$ diverges.

Proof.

(i) It can be seen from the proof process of lemma 1, when $0 < a < b+1$, $\int_0^1 \frac{\sin^b x}{x^a} dx$ converges. With the convergence of $\int_1^{+\infty} \frac{1}{x^a} dx$ and comparative discrimination, it's easy to know $\int_1^{+\infty} \frac{\sin^b x}{x^a} dx$ absolutely converges. So the conclusion is valid.

(ii) For $0 < a \le 1$, if $b = 2n-1$, because of the convergence of $\int_0^{+\infty} \frac{\sin x}{x^a} dx$ and equations in passage[4], $I$ converges, if $b = 2n$, by

$$\int_0^{+\infty} \frac{\sin^{2n} x}{x^a} dx = \int_0^{\frac{1}{2(n-k)}} \frac{\sin^{2n} x}{x^a} dx + \int_{\frac{1}{2(n-k)}}^{+\infty} \frac{\sin^{2n} x}{x^a} dx,$$

and

$$\frac{\sin^{2n} x}{x^a} \sim x^{2n-a} \quad (x \to 0^+),$$

so the first part is definite integral, the latter part of the integral involves $\int \frac{\cos 2(n-k)x}{x^a} dx$.

let $t = 2(n-k)x$, then $\int_{\frac{1}{2(n-k)}}^{+\infty} \frac{\cos 2(n-k)x}{x^a} dx = 2^{a-1}(n-k)^{a-1} \int_1^{+\infty} \frac{\cos x}{x^a} dx$, we can know the convergence from the Dirichlet discriminant of infinite integral, and notice that $\frac{1}{2^{2n}} C_{2n}^n \int_1^{+\infty} \frac{1}{x^a} dx$ diverges, so that $I$ diverges using equations listed before.

## III  Calculation formula of two kinds of generalized integral

In the following, we will give a lemma and two propositions to derive the calculation formulas of two kinds of generalized integrals related to $I$.

Lemma 2 (Euler's partial summation formula). Let $\{a_n\}$ and $\{b_n\}$ be two real number columns, $A_n = \sum_{k=1}^n a_k$, then

$$\sum_{k=1}^{n} a_k b_k = A_n b_n - \sum_{k=1}^{n-1} A_k (b_{k+1} - b_k).$$

Proposition 1 If the generalized integral, $I = \int_{0}^{+\infty} \frac{\sin^b x}{x^{a+1}} dx \ (a, b \in N^+)$, converges, then its value is

$$I = \frac{\left| \sum_{k=1}^{b} k^a t_k C_b^{\frac{b-k}{2}} d_k \right|}{2^{b-1} \Gamma(a+1)},$$

and $t_k = \dfrac{\left[1 - (-1)^{b-k}\right](-1)^{\left[\frac{k}{2}\right]}}{2}$, $d_k = \begin{vmatrix} \dfrac{\pi}{2} & \ln k \\ \cos\dfrac{a+b}{2}\pi & \sin\dfrac{a+b}{2}\pi \end{vmatrix}$.

Proof.

Because

$$I = \int_{0}^{+\infty} \frac{\sin^b x}{x^{a+1}} dx = -\left. \frac{\sin^b x}{ax^a} \right|_{0}^{+\infty} + \int_{0}^{+\infty} \frac{(\sin^b x)'}{ax^a} dx.$$

Substitute the Taylor expansion of sine function at 0: $\sin x = x + o(x)$, for $0 \le k \le a-1$,

$$\lim_{x \to 0} \frac{\sin^{b-k} x}{x^{a-k}} = \lim_{x \to 0} x^{b-a} + o(x^{b-a}) = 0 = \lim_{x \to \infty} \frac{\sin^{b-k} x}{x^{a-k}}.$$

Therefore, it can be seen from the integral method of parts and the inductive method of Mathematics that $I = \dfrac{1}{a!} \int_{0}^{+\infty} \dfrac{(\sin^b x)^{(a)}}{x} dx$.

Consider the Fourier series expansion of the numerator: $\sin^b x = \sum_{k=0}^{b} (A_k \cos kx + B_k \sin kx)$,

and from mathematical induction

$$(\sin^b x)^{(a)} = \sum_{k=0}^{b} k^a \left[ A_k \cos\left(kx + \frac{a}{2}\pi\right) + B_k \sin\left(kx + \frac{a}{2}\pi\right) \right].$$

Use trigonometric function formula,

$$I = \frac{\pi}{2a!}\sum_{k=1}^{b}k^{a}\left(B_{k}\cos\frac{a}{2}\pi - A_{k}\sin\frac{a}{2}\pi\right)$$

$$+\frac{1}{a!}\sum_{k=1}^{b}k^{a}\int_{0}^{+\infty}\frac{\left(A_{k}\cos\frac{a}{2}\pi + B_{k}\sin\frac{a}{2}\pi\right)\cos kx}{x}dx.$$

Let $C_{k} = A_{k}\cos\frac{a}{2}\pi + B_{k}\sin\frac{a}{2}\pi$, $I' = \int_{0}^{+\infty}\frac{\sum_{k=1}^{b}k^{a}C_{k}\cos(kx)}{x}dx$. It can be seen from the

Frullani integral that $\int_{0}^{+\infty}\frac{\cos(\alpha x)-\cos(\beta x)}{x}dx = \ln\frac{\beta}{\alpha}$.

So

$$I' = \int_{0}^{+\infty}\frac{1^{a}C_{1}(\cos x - \cos 2x) + (1^{a}C_{1} + 2^{a}C_{2})(\cos 2x - \cos 3x) + \cdots + \sum_{k=1}^{b}k^{a}C_{k}\cos bx}{x}dx$$

$$= \int_{0}^{+\infty}\frac{\sum_{k=1}^{b-1}k^{a}C_{k}[\cos kx - \cos(k+1)x] + \sum_{k=1}^{b}k^{a}C_{k}\cos bx}{x}dx$$

Easy to know exchange integral and sum order, for $\forall n \leq b-1$,

$$\int_{0}^{+\infty}\frac{\sum_{k=1}^{n}k^{a}C_{k}[\cos nx - \cos(n+1)x]}{x} = \sum_{k=1}^{n}k^{a}C_{k}\ln\frac{n+1}{n}.$$

If $C_{b} \neq 0$, because $\int_{0}^{+\infty}\frac{\cos bx}{x}dx$ diverges, it must satisfy that: $\sum_{k=1}^{b}k^{a}C_{k} = 0$.

Else if $C_{b} = 0$, the same way to get $\sum_{k=1}^{b-1}k^{a}C_{k} = 0$, we add $b^{a}C_{k}$ to both sides of equations,

so there always is: $\sum_{k=1}^{b}k^{a}C_{k} = 0$.

Thus, it is not difficult to get it from lemma 2

$$I' = \sum_{k=1}^{b}k^{a}C_{k}\ln b - \sum_{k=1}^{b}k^{a}C_{k}\ln k = -\sum_{k=1}^{b}k^{a}C_{k}\ln k.$$

Notice that $C_{k} = A_{k}\cos\frac{a}{2}\pi + B_{k}\sin\frac{a}{2}\pi$:

$$I = \sum_{k=1}^{b}\frac{k^{a}}{a!}\left[\frac{\pi}{2}\left(B_{k}\cos\frac{a}{2}\pi - A_{k}\sin\frac{a}{2}\pi\right) - \ln k\left(A_{k}\cos\frac{a}{2}\pi + B_{k}\sin\frac{a}{2}\pi\right)\right].$$

However, the power formula of trigonometric function sorted out in the previous paper finds the corresponding $A_{k}$, $B_{k}$, by comparing the two analysis results. After changing the element, it should meet the one-to-one correspondence, then there must be that $b$, $k$ both are odd or even numbers, so it can be written uniformly as

$$I = \frac{\sum_{k=1}^{b} k^a \frac{1-(-1)^{b-k}}{2} C_b^{\frac{b-k}{2}} (-1)^{\left[\frac{k}{2}\right]+1} d_k}{2^{b-1} a!},$$

(If parities of $b$, $k$ are different, $\frac{1-(-1)^{b-k}}{2}=0$ will make the equality before and after the change corresponds.) and $d_k = \begin{vmatrix} \frac{\pi}{2} & \ln k \\ \cos\frac{a+b}{2}\pi & \sin\frac{a+b}{2}\pi \end{vmatrix}$. Let $t_k = \frac{\left[1-(-1)^{b-k}\right](-1)^{\left[\frac{k}{2}\right]}}{2}$,

note the periodicity of trigonometric function, and use absolute value reduction to get

$$\int_0^{+\infty} \frac{\sin^b x}{x^{a+1}} dx = \frac{\left|\sum_{k=1}^{b} k^a t_k C_b^{\frac{b-k}{2}} d_k\right|}{2^{b-1} \Gamma(a+1)}.$$

Proposition 2. If the generalized integral,

$$I = \int_0^{+\infty} \frac{\sin^b x}{x^{a+\rho}} dx \, (a \in N, b \in N^+, \rho \in (0,1)),$$ converges, then its value is

$$I = \frac{\pi}{2^{b-1} \Gamma(a+\rho)} \left| \sum_{k=1}^{b} k^{a+\rho-1} t_k C_b^{\frac{b-k}{2}} \frac{\sin\frac{a+b+\rho}{2}\pi}{\sin \rho\pi} \right|,$$

and $t_k = \frac{\left[1-(-1)^{b-k}\right](-1)^{\left[\frac{k}{2}\right]}}{2}$.

Proof.

From proposition 1,

$$I = \int_0^{+\infty} \frac{\sin^b x}{x^{a+\rho}} dx = -\frac{\sin^b x}{(a+\rho-1) x^{a+\rho-1}} \Big|_0^{+\infty} + \int_0^{+\infty} \frac{(\sin^b x)'}{(a+\rho-1) x^{a+\rho-1}} dx,$$

for $k \leq b$,

$$\lim_{x \to 0} \frac{\sin^{b-k} x}{x^{a+\rho-1-k}} = \lim_{x \to 0} x^{b-(a+\rho-1)} + o(x^{b-(a+\rho-1)}) = 0 = \lim_{x \to +\infty} \frac{\sin^{b-k} x}{x^{a+\rho-1-k}}.$$

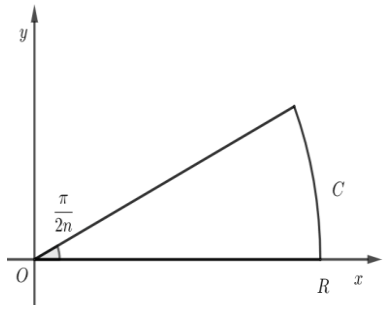

So $I = \dfrac{\Gamma(\rho)}{\Gamma(a+\rho)} \int_0^{+\infty} \dfrac{(\sin^b x)^{(a)}}{x^\rho} dx$.

Consider proving generalized Fresnel integral in [5] ($n>1$),

$$I_1 = \int_0^{+\infty} \sin x^n dx = \Gamma\left(1+\dfrac{1}{n}\right)\sin\dfrac{\pi}{2n}$$

$$I_2 = \int_0^{+\infty} \cos x^n dx = \Gamma\left(1+\dfrac{1}{n}\right)\cos\dfrac{\pi}{2n}.$$

Proof.

Take the sector with the center angle of $\dfrac{\pi}{2n}$ and the radius of infinity as the contour integral (as shown in the figure above).

$$\int_0^\infty e^{ix^n} dx = \int_0^\infty \cos x^n dx + i\sin\int_0^\infty \sin x^n dx = I_2 + iI_1.$$

Note that there is no singularity in the enclosure and we have $dz = e^{i\frac{\pi}{2n}} dR$.

So

$$\oint_C e^{iz^n} dz = \int_0^R e^{ix^n} dx + \int_{C_R} e^{iz^n} dz + \int_R^0 e^{iz^n} dz = 0,$$

$$z^n = R^n e^{i\frac{\pi}{2}} = R\left(\cos\dfrac{\pi}{2} + i\sin\dfrac{\pi}{2}\right) = iR^n,$$

$$\lim_{R\to\infty}\int_R^0 e^{iz^n} dz = -\int_0^\infty e^{i\cdot iR^n} \cdot e^{\frac{i\pi}{2n}} dR = -e^{\frac{i\pi}{2n}}\int_0^\infty e^{-R^n} dR.$$

Easy to know that $\int_{C_R} e^{iz^n} dz = 0$. With $\int_0^{+\infty}\left(\dfrac{x^{\frac{1}{n}}}{e^x}\right)' dx = 0$, so $\int_0^{+\infty} e^{-x^n} dx = \Gamma\left(1+\dfrac{1}{n}\right)$. Finally, the real part and the virtual part are taken respectively, and the result is obtained by the transformation of the generalized Fresnel integral,

$$I_\alpha = \int_0^{+\infty} \dfrac{\sin kx}{x^n} dx = \dfrac{\pi}{2}\dfrac{k^{n-1}}{\Gamma(n)\sin\dfrac{n}{2}\pi}$$

$$I_\beta = \int_0^{+\infty} \dfrac{\cos kx}{x^n} dx = \dfrac{\pi}{2}\dfrac{k^{n-1}}{\Gamma(n)\cos\dfrac{n}{2}\pi}.$$

So

$$I = \dfrac{\Gamma(\rho)}{\Gamma(a+\rho)}\sum_{k=1}^b k^a \int_0^{+\infty} \dfrac{\left(B_k\cos\dfrac{a}{2}\pi - A_k\sin\dfrac{a}{2}\pi\right)\sin kx + \left(A_k\cos\dfrac{a}{2}\pi + B_k\sin\dfrac{a}{2}\pi\right)\cos kx}{x^\rho} dx.$$

Then

$$I = \frac{\pi}{2\Gamma(a+\rho)} \sum_{k=1}^{b} k^{a+\rho-1} \left( \frac{B_k \cos\frac{a}{2}\pi - A_k \sin\frac{a}{2}\pi}{\sin\frac{\rho}{2}\pi} + \frac{A_k \cos\frac{a}{2}\pi + B_k \sin\frac{a}{2}\pi}{\cos\frac{\rho}{2}\pi} \right).$$

Again, we find the corresponding $A_k$ and $B_k$ from the power formula and simplify them, then we get when $a \in N, b \in N^+, \rho \in (0,1)$,

$$I = \frac{\pi}{2^{b-1}\Gamma(a+\rho)} \left| \sum_{k=1}^{b} k^{a+\rho-1} t_k C_b^{\frac{b-k}{2}} \frac{\sin\frac{a+b+\rho}{2}\pi}{\sin\rho\pi} \right|.$$

## IV  Unified calculation formula

How to unify the two situations into one expression? We hope that the formula will not distinguish the integer and decimal parts of the power of the denominator. It is easy to find that $\int_0^{+\infty} \frac{\cos kx}{x^\rho} dx$ diverges when $\rho = 1$, but converges when $\rho \in (0,1)$. Different convergences and divergences lead to different expressions. The unified expression of limit is considered. Notice that,

$$\lim_{\rho \to 1} \frac{\pi}{2} \frac{k^{\rho-1}-1}{\cos\frac{\rho}{2}\pi} = -\ln k$$

$$\lim_{\rho \to 1} \frac{\pi}{2} \frac{k^\rho - 1}{\sin\frac{\rho}{2}\pi} = \ln k$$

Let $C_k = A_k \cos\frac{a}{2}\pi + B_k \sin\frac{a}{2}\pi$, $D_k = B_k \cos\frac{a}{2}\pi - A_k \sin\frac{a}{2}\pi$, then:

$$\sum_{k=1}^{b} k^{a+\rho-1} \left( \frac{A_k \cos\frac{a}{2}\pi + B_k \sin\frac{a}{2}\pi}{\cos\frac{\rho}{2}\pi} + \frac{B_k \cos\frac{a}{2}\pi - A_k \sin\frac{a}{2}\pi}{\sin\frac{\rho}{2}\pi} \right).$$

$$= \sum_{k=1}^{b} k^a \left( C_k \frac{k^{\rho-1}-1}{\cos\frac{\rho}{2}\pi} + D_k \frac{k^{\rho-1}-k^{-1}}{\sin\frac{\rho}{2}\pi} + \frac{C_k}{\cos\frac{\rho}{2}\pi} + \frac{D_k k^{-1}}{\sin\frac{\rho}{2}\pi} \right)$$

The summation of the first two terms must be convergent which can be known from before. When $\rho \in (0,1)$, the sum of the last two terms is also convergent, when $\rho \to 1$, it can be transformed into the first kind of integral, if the integral converges, $\sum_{k=1}^{b} k^a C_k = 0$, which is consistent with previous conclusions, when $\rho \to 0$, if it converges, $\sum_{k=1}^{b} k^{a-1} D_k = 0$, after substitution, it is equivalent to the first kind of commutation.

To sum up, this paper first gives the convergence results of the generalized integral

$$I = \int_0^{+\infty} \frac{\sin^b x}{x^a} dx \text{ for } a \in R^+, b \in N^+,$$ then calculates the first kind of integral whose denominator power is integer and the second kind of integral whose denominator power is non integer in turn, and then uses the limit and groups it, so that the first kind of integral is actually the limit of the second kind of integral, and unifies the two, and finally gets it Unified calculation formula that, for $b+1 > a, a \in R^+, b \in N^+$ and $t_k = \dfrac{\left[1+(-1)^{b-k}\right](-1)^{\left[\frac{k}{2}\right]}}{2}$,

$$\int_0^{+\infty} \frac{\sin^b x}{x^a} dx = \frac{\pi}{2^{b-1}\Gamma(a)} \lim_{\alpha \to a} \left| \sum_{k=1}^{b} k^{\alpha-1} t_k C_b^{\frac{b-k}{2}} \frac{\sin\frac{\alpha+b}{2}\pi}{\sin \alpha\pi} \right|.$$

# [References]